\documentclass[11pt]{amsart}
\usepackage{amsfonts,amsmath}
\usepackage[dvips]{graphicx}
\newdimen\margin   
\def\COMMENT#1{}

\def\endproof{\noproof\bigskip}
\def\proof{\removelastskip\penalty55\medskip\noindent{\bf Proof. }}

\def\noproof{{\unskip\nobreak\hfill\penalty50\hskip2em\hbox{}\nobreak\hfill%
       $\square$\parfillskip=0pt\finalhyphendemerits=0\par}\goodbreak}

\newcommand{\eps}{\varepsilon}

\newcommand{\prob}{\mathbb{P}}
\newcommand{\ex}{\mathbb{E}}
\newcommand{\Pa}{{\mathcal P}}
\newcommand{\U}{{\mathcal U}}
\newcommand{\B}{{\mathcal B}}
\newcommand{\es}{{\mathcal S}}
\newcommand{\E}{{\mathcal E}}
\newcommand{\W}{{\mathcal W}}
\newcommand{\F}{{\mathcal F}}
\newcommand{\M}{{\mathcal M}}
\newcommand{\C}{{\mathcal C}}
\newcommand{\D}{{\mathcal D}}
\newcommand{\X}{{\mathcal X}}
\newcommand{\cG}{{\mathcal G}}

\newcommand{\Aux}{{\mathcal G}}
\newcommand{\aux}{{\mathcal G}}
\newcommand{\AuxU}{{\mathcal A}}

\newtheorem{firsttheorem}{Proposition}

\newtheorem{theorem}[firsttheorem]{Theorem}
\newtheorem{lemma}[firsttheorem]{Lemma}
\newtheorem{corollary}[firsttheorem]{Corollary}

\newtheorem{problem}[firsttheorem]{Problem}

\begin{document}
\title{Minors in random regular graphs}
\author{Nikolaos Fountoulakis, Daniela K\"uhn and Deryk Osthus}
\date{\today}
\begin{abstract} We show that there is a constant $c$ so that
for fixed $r \ge 3$ a.a.s.~an $r$-regular graph on $n$ vertices contains a complete
graph on $c \sqrt{n}$ vertices as a minor. This confirms a conjecture of Markstr\"om~\cite{Mark}.
Since any minor of an $r$-regular graph on $n$ vertices
has at most $rn/2$ edges, our bound is clearly best possible up to
the value of the constant $c$. As a corollary, we also obtain the likely order of magnitude
of the largest complete minor in a random graph $G_{n,p}$ during the phase transition
(i.e.~when $pn \to 1$).
\end{abstract}
\maketitle
\section{Introduction}
We say that a graph $G$ contains a complete graph on $k$ vertices (denoted by $K_k$)
as a minor if we can obtain a copy of~$K_k$ after a series of contractions of the
edges and deletions of vertices or edges of~$G$. We write $K_k \prec G$ in this case. Equivalently,
$G$ has a $K_k$~minor if there are~$k$ pairwise disjoint non-empty subsets of~$V(G)$
(which we call \emph{branch sets})
such that each of them is connected and any two of them are joined by an edge.
The \emph{contraction clique number}~ccl$(G)$ of~$G$ is the largest integer~$k$ such
that~$G$ has a $K_k$~minor.

Originally, the study of the order of the largest complete minor in a random graph
was motivated by Hadwiger's conjecture which states that ccl$(G) \geq \chi(G)$ for any graph~$G$.
Bollob\'as, Erd\H{o}s and Catlin~\cite{BCE} showed that the proportion of graphs on~$n$ vertices that
satisfy Hadwiger's conjecture tends to~$1$ as~$n$ tends to infinity. For this, they determined
the likely value of ccl$(G_{n,p})$ for the random graph~$G_{n,p}$ with constant edge probability~$p$
and compared this with known results on $\chi(G_{n,p})$.
Krivelevich and Sudakov~\cite{KS} investigated ccl$(G)$ for expanding graphs~$G$ and derived
the order of magnitude of {\rm ccl}$(G_{n,p})$ from their results when~$p$ is a polynomial in~$n$.
In~\cite{FKO}, we extended these results to any~$p$ with $pn \ge c$ for some constant $c>1$,
which answered a question from~\cite{KS}.
In particular, we showed that if $pn=c$ for some fixed $c>1$ then a.a.s.
\begin{equation} \label{sqrt}
{\rm ccl}(G_{n,p})=\Theta(\sqrt{n}).
\end{equation}
The upper bound is immediate, as for such~$p$ a.a.s.~the random graph
$G_{n,p}$ has $\Theta(n)$ edges and no minor of a graph~$G$ can contain more edges than~$G$
itself. Here we write that an event regarding a graph
on~$n$ vertices holds a.a.s.~if the probability of this event tends to~$1$ as~$n$ tends to infinity.

Markstr\"om~\cite{Mark} had earlier conjectured a similar phenomenon as in~(\ref{sqrt})
for the case of random regular graphs.
For any $r\geq 3$ and $n\geq 4$ such that~$rn$ is even, we denote by $G(n,r)$ a graph
chosen uniformly at random from the set of $r$-regular simple graphs on~$n$ vertices.
Throughout, we consider the case where~$r$ is fixed.
The number of edges of~$G(n,r)$ is $rn/2$ and so the same
argument as above shows that ccl$(G(n,r)) \le 2\sqrt{rn}$.
However the lower bound in~(\ref{sqrt}) does not imply that a random $r$-regular graph
satisfies {\rm ccl}$(G(n,r))=\Omega(\sqrt{n})$ as the asymptotic structure of $G(n,r)$ is quite
different from that of $G_{n,{r/n}}$ (see for example~\cite{Worm} or Chapter 9 in~\cite{JLR}).
Markstr\"om~\cite{Mark} proved that $G(n,3)$ a.a.s.~contains a complete minor of order~$k$,
for any integer $k\geq 3$ and conjectured that $G(n,3)$ contains a complete minor of order~$\Omega(\sqrt{n})$.
In this paper, we verify this conjecture for any $r\geq 3$:
\begin{theorem} \label{MainThm}
There exists an absolute constant $c>0$ such that for every fixed $r\geq 3$ a.a.s.
$c\sqrt{n}\leq {\rm ccl}(G(n,r)) \leq 2\sqrt{rn}$.
\end{theorem}
This result can be combined with results of {\L}uczak~\cite{Luc91} to determine the
likely order of magnitude of {\rm ccl}$(G_{n,p})$ during the phase transition, i.e.~when $pn \to 1$
(see Section~\ref{secphase} for details).
\begin{corollary} \label{phase}
There exists an absolute constant $c>0$ such that whenever
$np=1+\lambda n^{-1/3}$, where $\lambda =\lambda (n)\to \infty$ but $\lambda =o(n^{1/3})$,
then a.a.s. $c\lambda^{3/2} \leq {\rm ccl}(G_{n,p}) \leq 4 \lambda^{3/2}$.
\end{corollary}
{\L}uczak, Pittel and Wierman~\cite{LPW} previously showed that a.a.s.~ccl$(G_{n,p})$ is unbounded for~$p$
as in Corollary~\ref{phase}.
For smaller $p$ (i.e.~when $np \le 1+\lambda n^{-1/3}$ for some constant $\lambda$)
they showed that~${\rm ccl}(G_{n,p})$ is bounded in probability, i.e. for every $\delta >0$ there exists $C=C(\delta)$ such that
$\prob ({\rm ccl}(G_{n,p})> C)< \delta$.
As described earlier, values of $p$ which are larger than those allowed for in Corollary~\ref{phase}
but bounded away from~$1$ are covered in~\cite{FKO}.
So altogether, all these results
determine the likely order of magnitude of
ccl$(G_{n,p})$ for any~$p$ which is bounded away from~$1$.

The results in~\cite{LPW} were proved
in connection with the following result on the limiting probability $g(\lambda)$ that
$G_{n,p}$ is planar in the above range. The authors proved that if~$\lambda$ is bounded,
then~$g(\lambda)$ is bounded away from~$0$ and~$1$.
If $\lambda \to - \infty$, then $g(\lambda) \to 1$, whereas
if $\lambda  \to \infty$, then $g(\lambda) \to 0$.

The result in~(\ref{sqrt}) and Theorem~\ref{MainThm} raise the question of
whether one can extend these results to other (not necessarily random) graphs.
A natural class to consider are expanding graphs: A graph $G$ on $n$ vertices
is an \emph{$(\alpha,t)$-expander} if any $X \subseteq V(G)$ with $|X| \le \alpha n/t$ satisfies
$|N(X)| \ge t|X|$, where $N(X)$ denotes the external neighbourhood of~$X$.
\begin{problem} \label{problem}
Is there a constant $c>0$ such that for each~$r \ge 3$ every
$r$-regular $(1/3, 2)$-expander satisfies ${\rm ccl}(G)\ge c\sqrt{n}$?
\end{problem}
An answer to the problem would indicate whether expansion alone is sufficient
when trying to force a complete minor of the largest possible order in a sparse graph, or
whether other parameters are also relevant.
Krivelevich and Sudakov~\cite{KS} showed that we do have ccl$(G)\ge c\sqrt{n/\log n}$ if $r \ge 10$.
(They also considered the case when $r$ is not bounded but grows with $n$.)
As observed in~\cite{KS}, this bound can also be deduced from a result of Plotkin, Rao and Smith~\cite{PRS}
on separators in graphs without a large complete minor.
Kleinberg and Rubinfeld~\cite{KR} also considered the same problem but with a weaker definition of expansion.

One can ask similar questions as above for topological minors.
Topological minors in random graphs were investigated in~\cite{BC,BCE,AKS}.
An analogue of Problem~\ref{problem} would be to ask for which values of $\alpha$, $t$, $r$
an $r$-regular $(\alpha,t)$-expander on $n$ vertices contains a subdivision of a $K_{r+1}$.
We expect that this might not be difficult to prove for fixed $r$ but harder if $r$ is no
longer very small compared to $n$.

\section{Proof of Corollary~\ref{phase}} \label{secphase}
The upper bound in Corollary~\ref{phase} will follow from basic facts about minors as well as
the structure of~$G_{n,p}$. Bollob\'as~\cite{Bol84} (see also~\cite{Bol} or~\cite{JLR}) proved
that a.a.s.~all the components of~$G_{n,p}$, except from the largest one, are either trees or unicyclic.
Therefore none of them contains a $K_4$~minor. Let $L_1(G_{n,p})$
denote the largest component of~$G_{n,p}$.
Given a graph~$G$, we define its \emph{excess} as ${\rm exc}(G):=e(G)-|G|+1$. (${\rm exc}(G)$ is
also called the \emph{cyclomatic number} of~$G$.)
Observe that if~$H$ and~$G$ are connected graphs and $H \prec G$ then ${\rm exc}(H)\leq {\rm exc}(G)$.
Since ${\rm exc}(K_r)={r\choose 2} - r+1 \geq r^2/16$ for $r\geq 4$
this implies that if $K_r \prec L_1(G_{n,p})$ for some $r\geq 4$ then $r^2/16
\leq {\rm exc} (L_1(G_{n,p}))$, or equivalently
\begin{equation} \label{eq:excessSize}
r \leq 4 \sqrt{{\rm exc}(L_1 (G_{n,p}))}.
\end{equation}
{\L}uczak~\cite{Luc90} gave a tight estimate on ${\rm exc}(L_1 (G_{n,m}))$, where~$G_{n,m}$
is a random graph with~$n$ vertices and~$m$ edges (i.e.~$G_{n,m}$ is chosen uniformly at random
among all such graphs). He proved that if $m=n/ 2 + \bar{\lambda} n^{2/3}$, where
$\bar{\lambda} = \bar{\lambda} (n) \to \infty$ but $\bar{\lambda}=o(n^{1/3})$, then a.a.s.
${\rm exc}(L_1 (G_{n,m}))=(1+o(1))16 \bar{\lambda}^3/3$.
This trivially implies that if $m=n/2+\bar{\lambda}n^{2/3}+O(\sqrt{n})$ then
a.a.s. ${\rm exc}(L_1 (G_{n,m})) \leq 8 \bar{\lambda}^3$.
Together with the fact that ${n \choose 2}p = n/2 + \lambda n^{2/3}/2 + O(1)$
and Proposition~1.12 in~\cite{JLR} this implies that a.a.s. ${\rm exc}(L_1 (G_{n,p}))\le \lambda^3$.
But if the latter holds and $K_r \prec L_1 (G_{n,p})$ for some $r\geq 4$ then~(\ref{eq:excessSize})
gives $r\leq 4 \lambda^{3/2}$. Thus a.a.s. ${\rm ccl}(G_{n,p}) \leq 4 \lambda^{3/2}$.

For the lower bound in Corollary~\ref{phase} we will use the following result of {\L}uczak
which is contained in the proof of Theorem~$5^*$ in~\cite{Luc91}.%
    \COMMENT{The $C$ in the thm actually lies $L_1(G_{n,m})$. But the statement of the thm is nevertheless
ok.}
\begin{theorem}\label{thm:subdcub}
Suppose that $m=n/ 2 + \lambda n^{2/3}$, where $\lambda \to \infty$ and $\lambda =o(n^{1/3})$.
Then there is a procedure which in any given graph~$G$ with~$n$ vertices and~$m$ edges finds a subdivision
of a (possibly empty) $3$-regular graph~$C=C(G)$ such that a.a.s. $|C(G_{n,m})|=(32/3+o(1))\lambda^3$
and conditional on $|C(G_{n,m})|=s$ in this range the distribution of~$C(G_{n,m})$ is the same as~$G(s,3)$.
\end{theorem}
Loosely speaking, Theorem~\ref{thm:subdcub} implies that a.a.s.~$L_1(G_{n,m})$ contains a subdivision of a
random 3-regular graph~$G(s,3)$ where $s=(32/3+o(1))\lambda^3$.  %
Together with Theorem~\ref{MainThm} this implies that a.a.s.
$${\rm ccl}(G_{n,m}) \geq {\rm ccl}(C(G_{n,m}))\ge  c \lambda^{3/2}.
$$
Again Proposition~1.12 from~\cite{JLR} now yields the lower bound of Corollary~\ref{phase}.

\section{Sketch of proof of Theorem~\ref{MainThm}}\label{sec:sketch}

We will use a result of Janson~\cite{janson} which implies that it suffices to find a complete minor
in the union of a random Hamilton cycle and a random perfect matching.
We split the  Hamilton cycle into paths~$P_1$ and~$P_2$ of equal length.
We further split~$P_1$ into $k$ connected candidate branch sets $B_i$, where~$k$ is close to~$\sqrt{n}$.
Each of these candidate branch sets has size roughly~$\sqrt n$. We now split~$P_2$ into
sets~$\Pa_i$ of disjoint paths. The lengths of the paths in~$\Pa_i$ is roughly~$3^i$, whereas the number
of paths in~$\Pa_i$ is roughly~$n/9^i$. For each pair $(B,B')$ of candidate branch sets we aim to find
a path~$P$ in some~$\Pa_i$ such that both~$B$ and~$B'$ are joined to~$P$ by an edge of the
random perfect matching. We let~$\U_{i-1}$ denote the set of pairs of candidate branch sets for
which we were not able to find such a path~$P$ in~$\bigcup_{j<i}\Pa_j$.
We will show inductively that $|\U_i| \le |\U_{i-1}|/27$ (with sufficiently high probability).
By continuing this for $(\log_3 n)/6$ stages and discarding a few atypical
branch sets, we eventually obtain the desired minor.
This strategy is similar to that of~\cite{FKO}. However, the proof that it works is very different:
the argument in~\cite{FKO} was based on a greedy matching algorithm whose analysis
crucially relied on the independence of certain events.
In the current setting, this no longer works. So instead, in each stage we use
Hall's theorem to find a large matching in the bipartite auxiliary graph whose vertex
classes are~$\U_{i-1}$ and~$\Pa_i$ and where a pair $(B,B')$ is adjacent to a path~$P\in \Pa_i$
if~$P$ can be used to join~$B$ and~$B'$ as above.
(Actually, it turns out that we need to consider suitable subsets $\U_{i-1}' \subseteq \U_{i-1}$ and
$\es \subseteq \Pa_i$ for the argument to work.)
Though the number~$|\Pa_i|$ of paths decreases in each stage, the increasing path length means that
the average degree of a pair $(B,B')$ in this auxiliary graph remains large
(but bounded) in each stage and so we can indeed expect
to find a large matching. On the other hand, one can show that there might be a significant number
of pairs from~$\U_{i-1}$ which are isolated in the auxiliary graph. So we cannot hope to get away with
just a single stage.

\section{Models of random $r$-regular graphs}
The aim of this section is to show that it suffices to find our complete minor in the
union of a random Hamilton cycle and a random perfect matching.
To do this, let us first describe the configuration model which was introduced by
Bender and Canfield~\cite{BenCan} and independently by Bollob\'as~\cite{Bol1}.
For $n\geq 1$ let $V_n:=\{1,\ldots, n \}$. Also for
those $n$ for which $rn$ is even, we let $P:=V_n \times [r]$.
A \emph{configuration} is a perfect matching on~$P$.
If we project a configuration onto~$V_n$, then we obtain an $r$-regular
multigraph on~$V_n$. Let~$G^*(n,r)$ denote the random multigraph
that is the projection of a configuration on~$P$ which is chosen uniformly at random. It
can be shown (see e.g.~\cite[p.~236]{JLR}) that if we condition on~$G^*(n,r)$
being simple (i.e.~it does not have loops or multiple edges), then this is distributed
uniformly among the $r$-regular graphs on~$V_n$. In other words, $G^*(n,r)$ conditional
on being simple has the same distribution as~$G(n,r)$. We also let~$G'(n,r)$
denote a random multigraph whose distribution is that of~$G^*(n,r)$
conditional on having no loops.
We will use the above along with the following (see Corollary 9.7
in~\cite{JLR}):
\begin{equation}\label{eq:LimSimple}
\lim_{n\rightarrow \infty} \prob (G^*(n,r) \mbox{ is simple} ) >0.
\end{equation}
(Of course the above limit is taken over those~$n$ for which~$rn$ is
even.) Let~$A_n$ be a subset of the set of $r$-regular multigraphs on~$V_n$.
Altogether the above facts imply that
if $\prob(G'(n,r)\in A_n)\rightarrow 0$ as $n\rightarrow \infty$
then $ \prob (G(n,r)\in A_n)\rightarrow 0$. Indeed, suppose that
the former holds. Then
\begin{align} \label{eq:Models}
\prob (G(n,r) & \in A_n) = \prob(G^*(n,r)\in A_n \,|\, G^*(n,r) \text{ simple})
=\frac{\prob(G^*(n,r)\in A_n, \ G^*(n,r) \mbox{ simple})}
{\prob (G^*(n,r) \mbox{ simple})} \nonumber \\
&\leq \frac{\prob(G^*(n,r)\in A_n, \ G^*(n,r) \mbox{ has no loops})}
{\prob (G^*(n,r) \mbox{ has no loops})\prob (G^*(n,r) \mbox{ simple})} =
\frac{\prob (G'(n,r) \in A_n)}{\prob (G^*(n,r) \mbox{ simple})}
\stackrel{(\ref{eq:LimSimple})}{\rightarrow} 0.
\end{align}
This allows us to work with~$G'(n,r)$ instead of~$G(n,r)$ itself.

Let us first assume that $r=3$.
The reason for working with~$G'(n,3)$ is that we may think of
it as being the union of a random Hamilton cycle on~$V_n$
and a random perfect matching on~$V_n$. This is made
precise by the notion of \emph{contiguity}. If~$(\mu_n)$ and~$(\nu_n)$
are two sequences of probability measures such that for each~$n$,
$\mu_n$ and~$\nu_n$ are measures on the same measurable space~$\Omega_n$,
then we say that they are \emph{contiguous} if for
every sequence of measurable sets $(A_n)$ with $A_n \in \Omega_n$ we
have $\lim_{n\rightarrow \infty} \mu_n (A_n) =0$ if and only if
$\lim_{n\rightarrow \infty} \nu_n (A_n)=0$.
Now let $H(n)+G(n,1)$ denote the random multigraph on~$V_n$ that is
obtained from a Hamilton cycle on~$V_n$ chosen uniformly at random by adding a random perfect
matching on~$V_n$ chosen independently from the Hamilton cycle.
Janson~\cite{janson} (see also Theorem 9.30 in~\cite{JLR}) proved that $H(n)+G(n,1)$ is
contiguous to~$G'(n,3)$.%
    \COMMENT{Makes sense as both graphs are defined iff $n$ is even}
\begin{theorem}\label{thm:contig}
The random 3-regular multigraphs $H(n)+G(n,1)$ and $G'(n,3)$ are contiguous.
\end{theorem}
So instead of proving Theorem~\ref{MainThm} directly, it suffices to prove the following result.
\begin{theorem} \label{Thm_3Reg}
There exists an absolute constant $c'>0$ such that a.a.s.~the random multigraph
$H(n)+G(n,1)$ contains a complete minor of order at least $c'\sqrt{n}$.
\end{theorem}
Together with~(\ref{eq:Models}) and Theorem~\ref{thm:contig} this then implies the
lower bound of Theorem~\ref{MainThm} for $r=3$.
The lower bound for $r>3$ follows from Theorem~9.36 in~\cite{JLR} which states that
for each $s\ge 3$ an increasing property that holds a.a.s.~for $G(n,s)$ also holds
a.a.s.~for $G(n,s+1)$.

\section{Notation and large deviation inequalities}
\subsection{Notation}
Given a graph~$G$ and two disjoint sets~$A$ and $B$ of vertices, we say that
an edge of~$G$ is an \emph{$A$-$B$ edge} if it joins a vertex in~$A$ to
a vertex in~$B$. Given disjoint subgraphs~$H$ and $H'$ of~$G$, we define
$H$-$H'$ edges of~$G$ similarly.
Given $a,b\in\mathbb{R}$ we write $[a\pm b]$ for the interval $[a-b,a+b]$.
We will write $\ln ^2 n$ for $(\ln n)^2$. We omit floors and ceilings whenever this
does not affect the argument.

\subsection{A concentration inequality}
In this subsection, we will state a concentration inequality which
we will use several times during the proof of Theorem~\ref{Thm_3Reg}.
This is Theorem~7.4 in~\cite{McD}. We first describe the more general
setting to which this theorem applies.

Let~$W$ be a finite probability space that is also a metric space with
its metric denoted by~$d$. Suppose that $F_0,\ldots ,F_s$ is a sequence of
partitions of~$W$ such that~$F_{j+1}$ refines~$F_j$, $F_0$ is the partition
consisting of only one part (i.e.~$F_0=\{W\}$) and~$F_s$ is the partition
where each part is a single element of~$W$.
Suppose that whenever $A, B \in F_{j+1}$ and $C\in F_j$
are such that $A,B \subseteq C$, then there is a bijection $\phi: A
\rightarrow B$ such that $d(x,\phi (x))\leq c$. Now, let $w\in W$ be chosen uniformly at
random and let $f:W\rightarrow \mathbb{R}$ be a function on~$W$ satisfying
$|f(x)-f(y)|\leq d(x,y)$. Then for all $a>0$
\begin{equation}\label{ConcM}
\prob \left( |f(w)- \ex\left(f(w)\right)|> a \right) \leq
2 \exp\left({-2a^2\over s c^2
}\right).
\end{equation}
\subsection{The hypergeometric distribution}
Let~$Z$ be a non-empty finite set and $Z'\subseteq Z$. Assume that we sample a set~$Y$
uniformly at random among all subsets of~$Z$ having size~$y$. Recall that the size of
$Y\cap Z'$ is a random variable whose distribution is \emph{hypergeometric} and whose
expected value is $\lambda:=y|Z'|/|Z|$. We will often use the following concentration
inequality that follows e.g.~from Theorem~2.10 and Inequalities~(2.5) and~(2.6) in~\cite{JLR}:
\begin{equation} \label{eq:Hyp_Conc1}
\prob \left(||Y\cap Z'| - \lambda | \geq  a \right) \leq
2\exp \left(-{a^2\over 2(\lambda + a/3) } \right)
\end{equation}
for all $a\geq 0$.

\section{Proof of Theorem~\ref{Thm_3Reg}}
\subsection{Setup}
Let~$V_n$ be a set of~$n$ vertices.
We will expose the random multigraph $H(n)+G(n,1)$ on~$V_n$ in stages starting with the
Hamilton cycle~$H(n)$. We split~$H(n)$ into two paths $P_1, P_2$ of equal lengths each
having~$n/2$ vertices.%
     \COMMENT{ok since $G(n,1)$ is only defined if $n$ is even}
As described in Section~\ref{sec:sketch}, the (candidate) branch sets for our minor
will be subpaths of~$P_1$ and we will use the edges of the random perfect matching~$G(n,1)$
as well as subpaths of~$P_2$ to join them.
Let us now turn to~$G(n,1)$. So consider a perfect matching~$M^*$ on~$V_n$ chosen uniformly
at random. Our first aim is to estimate the number of $P_1$-$P_2$ edges of~$M^*$.

\begin{lemma}\label{Xrange} With probability $1-O( 1/\ln^2 n)$
the number of $P_1$-$P_2$ edges of~$M^*$ lies in the interval $\left[n/4 \pm \sqrt{n}\ln n \right]$.
\end{lemma}
\proof
This is a simple application of Chebyshev's inequality.
For each vertex $v \in V(P_1)$ set $X_v:=1$ if $M^*$ matches $v$ to a vertex of~$P_2$
and set $X_v:=0$ otherwise. Then $X:=\sum_{v \in V(P_1)} X_v$ is the number of $P_1$-$P_2$ edges of~$M^*$.
Note that for every $v$ we have $$ \prob (X_v=1) = \frac{n/2}{n-1}=1/2+O(1/n).$$
So $\ex X=(n/2) (1/2+O(1/n))=n/4+O(1)$.
Also, for distinct $v,w \in V(P_1)$ we have
$$\prob (X_v=1 \mid X_w=1) =\frac{n/2-1}{n-3} =1/2+O(1/n).$$ This implies that
\begin{align*}
\ex (X^2) & = \sum_{v \in V(P_1)} \prob (X_v=1) +\sum_{v \neq w \in V(P_1)} \prob ( X_v=X_w=1)\\
& =\ex X +\frac{n}{2} \left(\frac{n}{2} -1 \right) \left( \frac{1}{4} +O(1/n) \right)
= \frac{n^2}{16}+O(n)=( \ex X )^2 ( 1+O(1/n)).
\end{align*}
So Chebyshev's inequality implies that
$$\prob \left( |X -n/4| \ge \sqrt{n} \ln n \right)
\le  \prob \left( |X -\ex X| \ge (\sqrt{n} \ln n)/2 \right)
= \frac{O(1/n) (\ex X )^2}{n \ln^2 n} =O(1/\ln^2 n),
$$
as required.
\endproof

Fix a positive constant~$\eps$. Throughout the proof we will assume that~$\eps$ is
sufficiently small for our estimates to hold. (All conditions on~$\eps$ will involve only
absolute constants, i.e.~will be independent of~$n$.) Suppose that~$n$ is sufficiently large
compared to~$1/\eps$. Let~$k$ and~$t$ be integers such that
\begin{equation}\label{eq:defkt}
{k \choose 2} = \eps^4 n \ \ \text{and} \ \ t:=\frac{\sqrt{n}}{\eps}.
\end{equation}
So $k=(1+o(1)) \eps^2 \sqrt{2n}$.
Consider any $X_1\subseteq V(P_1)$ and $X_2 \subseteq V(P_2)$
such that $|X_1|=|X_2|\in [n/4 \pm \sqrt{n}\ln n]$. Let~$X'_1\subseteq X_1$
be the set of the first~$kt$ vertices on~$P_1$ in~$X_1$. Let~$X'_2\subseteq X_2$ be any subset
of size~$kt$. Let~$\X$ denote the event that~$X_1$ and~$X_2$ are the set
of endvertices of the $P_1$-$P_2$ edges in our random perfect matching~$M^*$ on~$V_n$.
Similarly, let~$\X'$ be the event that~$M^*$ matches~$X'_1$ to $X'_2$.
In what follows, we will condition on both $\X$ and~$\X'$.
All our probability bounds will hold regardless of what the sets $X_1,X_2,X'_1,X'_2$
actually are (provided that~$|X_1|=|X_2|$ is in the specified range).

Pick~$k$ consecutive disjoint subpaths $B_1, \ldots , B_k$ of $P_1$ such that
$|V(B_i)\cap X_1|=|V(B_i)\cap X'_1|=t$ for all $i=1,\dots,k$. The $B_i$'s will be called \emph{candidate
branch sets} and the vertices in $V(B_i)\cap X'_1$ will
be called the \emph{effective vertices of~$B_i$}.
We will show that a.a.s. there is a complete minor whose branch sets%
     \COMMENT{Strictly speaking this is not correct as the branch sets will be bigger}
are almost all the~$B_i$'s.
Set
\begin{equation} \label{eq:i0def}
i_0 := (\log_3 n)/6.
\end{equation}
Choose consecutive disjoint subpaths $Q_1,\ldots , Q_{i_0}$ of~$P_2$ such that
\begin{equation} \label{eq:Qi}                                  
|Q_i|_{\rm eff}:=|V(Q_i)\cap X'_2|= {|X'_2|\over 3^i}={kt\over 3^i}={(1+o(1))\sqrt{2}\eps n\over 3^i}.
\end{equation}
The vertices in $V(Q_i)\cap X'_2$ are the \emph{effective vertices of~$Q_i$}
and $|Q_i|_{\rm eff}$ is the \emph{effective length of~$Q_i$}.
We further divide each~$Q_i$ into a set~$\Pa_i$ of
consecutive disjoint subpaths, each of effective length
\begin{equation} \label{eq:lidef}
\ell_i  := 100 \cdot 3^{i-1}.
\end{equation}
(So each of these subpaths meets~$X'_2$ in precisely~$\ell_i$ vertices.)
Note that
\begin{equation}\label{eq:li0}
\ell_{i_0} \leq 100 \cdot 3^{i_0} = 100 n^{1/6}
\end{equation}
and
\begin{equation} \label{eq:U_isize}
|\Pa_i|= {|Q_i|_{\rm eff}\over \ell_i} = {kt\over 300\cdot 9^{i-1}}.
\end{equation}
Thus $|\Pa_{i_0}|=\Theta(n^{2/3})$.
The strategy of our proof is to expose the neighbours of the (effective) vertices from~$X_2'$ in our
random perfect matching~$M^*$ in~$i_0$ stages. More precisely, during
the $i$th stage we will expose the neighbours of the effective vertices in~$Q_i$
(for every $1\leq i \leq i_0$).
We will show that with high probability during each stage we can use the paths in~$\Pa_i$
to join a large proportion of all those pairs of candidate branch sets that are still unjoined
after the previous stages. More precisely, an unjoined pair $(B,B')$ of candidate branch sets can
be joined through $P\in \Pa_i$ if our random perfect matching~$M^*$ contains both a $B$-$P$ edge
and a $B'$-$P$ edge. In this case we will say that~$P$ \emph{can be used to join the pair $(B,B')$}.
Of course, if we use~$P$ to join $(B,B')$ then~$P$ cannot be used to join another
unjoined pair of candidate branch sets.

Let us make the above more precise. Given $1\leq i \leq i_0$, let $\U_{i-1}$ denote the set
of pairs of candidate branch sets that are still unjoined after the $(i-1)$th stage. So~$\mathcal{U}_0$
is the set of all pairs of candidate branch sets. Note that
\begin{equation} \label{eq:U_0Value}
U_0:=|\mathcal{U}_0|={k \choose 2}=\eps^4 n.
\end{equation}
We will show that with high probability during the $i$th stage
we can join $26 |\U_{i-1}|/ 27$ pairs in~$\U_{i-1}$
using the paths belonging to~$\Pa_i$. So inductively we will prove that with high probability
\begin{equation}\label{eq:InductiveStep}
U_{i}:=|\U_{i}|={U_0 \over 27^{i}}=\frac{\eps^4 n}{27^i}.
\end{equation}
Suppose that~(\ref{eq:InductiveStep}) holds for all $j<i$ and that we now wish to analyze the~$i$th stage.
It will turn out that the pairs in~$\U_{i-1}$ which contain candidate branch sets
lying in too many other pairs from~$\U_{i-1}$ create problems. So we will ignore these
pairs. More precisely, let~$\B_{i-1}$ be the set of all those candidate branch sets that belong to
more than%
     \COMMENT{$\Delta_{i-1}=\Omega(\sqrt{n}/(2\cdot 9)^{i_0})=\Omega(n^{1/6}/3^{i_0\log_3 2})
=\Omega(n^{1/6}/n^{(\log_3 2)/6})$ is a large number}
\begin{equation} \label{eq:def_of_delta}
\Delta_{i-1}:= {(3/2)^{i-1} U_{i-1}\over  \eps^{1/8} k}=\frac{U_0}{\eps^{1/8}(2\cdot 9)^{i-1}k}
\end{equation}
pairs in~$\U_{i-1}$. Note that since
$|\B_{i-1}|\Delta_{i-1}\leq 2 U_{i-1}$ we have
\begin{equation}\label{eq:Bi}
|\B_{i-1}|\leq {2U_{i-1}\over \Delta_{i-1}}\le
{2\eps^{1/8}k \over (3/2)^{i-1}}.
\end{equation}
Let~$\U^*_{i-1}$ be the set of all those pairs in~$\U_{i-1}$ having at least one branch set
in~$\B_{i-1}$. Call these pairs \emph{bad}. If $|\U^*_{i-1}|\ge 26 U_{i-1}/27$, delete
precisely $26 U_{i-1}/27$ bad pairs from~$\U_{i-1}$ to obtain~$\U_i$.
If $|\U^*_{i-1}|< 26 U_{i-1}/27$ we let $\U'_{i-1}:=\U_{i-1}\setminus \U^*_{i-1}$.
We will show that during the $i$th stage with high probability we can join all but
$U_{i-1}/27$ pairs in~$\U'_{i-1}$. We let~$\U_i$ be the
set of the remaining unjoined pairs in $\U'_{i-1}$. Thus in both cases~$U_i=|\U_i|$
satisfies~$(\ref{eq:InductiveStep})$ with high probability.

After the end of stage~$i_0$ will delete all the candidate branch sets
in $\B_0\cup\dots\cup B_{i_0-1}$ (see Section~\ref{sec:final}). The number of these
candidate branch sets is
\begin{equation}\label{eq:TotalB_i}
\sum_{i=1}^{i_0} |\B_{i-1}| \stackrel{(\ref{eq:Bi})}{\leq}
 \sum_{i\geq 1}{2\eps^{1/8}k \over (3/2)^{i-1}} = 6\eps^{1/8} k.
\end{equation}

\subsection{Bounds on the number of effective vertices still available}
We will now estimate the number of all those effective vertices in each candidate branch
set that are joined to a path in~$\Pa_1\cup\dots\cup \Pa_{i-1}$, i.e.~that are matched after
the first $i-1$ stages. The total number of effective vertices in the candidate branch sets
that are matched after the first~$i-1$ stages is%
     \COMMENT{We have $\sum_{s=1}^{i-1} 3^{-s} = {1\over 3} \sum_{s=0}^{i-2}3^{-s}
={1\over 3} {1-3^{-i+1}\over 1-1/3} = {1\over 3} {1-3^{-i+1}\over 2/3}=
{1\over 2} \left(1-3^{-i+1}\right)$.}
\begin{equation}\label{eq:available}
\sum_{j\le i-1}|Q_j|_{\rm eff}
\stackrel{(\ref{eq:Qi})}{=} kt\sum_{j=1}^{i-1}{1\over 3^j}={kt \over 2}
~\left(1-3^{-(i-1)} \right)=:x_{i-1}.
\end{equation}
Each $x_{i-1}$-subset of the union~$X'_1$ of all the effective vertices in the candidate
branch sets is equally likely to be the set of these matched vertices.
Thus for every candidate branch set~$B$ the distribution of the number ${\rm eff}'_{i}(B)$
of all those effective vertices in~$B$ which are matched to (effective vertices on)
paths in $\Pa_1\cup \dots\cup \Pa_{i-1}$ is hypergeometric. Since in total~$B$ contains~$t$ effective
vertices and $|X'_1|=kt$ we can now use~(\ref{eq:Hyp_Conc1}) to see that
$$
\prob \left(|{\rm eff}'_{i}(B)-x_{i-1} t/kt |\geq n^{1/4} \ln n \ | \
\X,\X' \right)\le \exp\left(-\Omega(\ln^2 n)\right).$$
Thus,
\begin{equation*} 
{\rm eff}'_{i}(B) \in \left[ {x_{i-1}\over k} \pm n^{1/4}\ln n \right] \subseteq
\left[{t\over 2}(1-3^{-(i-1)}) \pm n^{1/3} \right],
\end{equation*}
with (conditional) probability $1-\exp\left(-\Omega(\ln^2 n)\right)$.
Now, let ${\rm eff}_i(B) := t - {\rm eff}'_i(B)$ be the number of all those effective vertices in~$B$
that are still unmatched after the first~$i-1$ stages and let ${\rm Eff}_i(B)$ denote the set of
all those effective vertices. Thus with (conditional) probability
$1-k\exp\left(-\Omega(\ln^2 n)\right)=1-\exp\left(-\Omega(\ln^2 n)\right)$
we have
\begin{equation} \label{eq:e_2j}
{\rm eff}_i(B) \in \left[ {t\over 2}(1+3^{-(i-1)}) \pm n^{1/3}\right]
\end{equation}
for all candidate branch sets~$B$.

Let~$M^*_{i-1}$ be any matching which matches the set ${\rm Eff}(Q_1)\cup\dots\cup {\rm Eff}(Q_{i-1})$
of effective vertices on the paths $Q_1,\dots,Q_{i-1}$ (equivalently the set
of effective vertices on the paths in $\Pa_1\cup\dots\cup \Pa_{i-1}$) into the
set of effective vertices in the candidate branch sets. Suppose that~$M^*_{i-1}$
is the submatching of our random matching~$M^*$ exposed after the first~$i-1$ stages.
Then~$M^*_{i-1}$ determines ${\rm Eff}_i(B)$ for every candidate branch set~$B$.
Moreover, by considering a fixed ordering of all the pairs in~$\U_0$, we may assume
that~$M^*_{i-1}$ also determines~$\U_{i-1}$. Call~$M^*_{i-1}$ \emph{good} if~(\ref{eq:e_2j})
holds for all candidate branch sets~$B$ and if~(\ref{eq:InductiveStep}) holds for~$i-1$.
Consider any good~$M^*_{i-1}$ and let~$\M^*_{i-1}$ denote the event that~$M^*_{i-1}$
is the submatching of our random matching~$M^*$ exposed after the first~$i-1$ stages.
From now on we will condition on~$\X$, $\X'$ and $\M^*_{i-1}$ and we let~$\prob_i(\cdot)$
denote the corresponding conditional probability measure that arises from choosing
a random matching from the set~${\rm Eff}(Q_i)$ of effective vertices on~$Q_i$ into the set
$\bigcup_{j=1}^k {{\rm Eff}_i(B_j)}$ of all those effective vertices in the
candidate branch sets which are not already endvertices of edges in~$M^*_{i-1}$
(i.e.~into the set of all those effective vertices
in the candidate branch sets that are still unmatched after the first~$i-1$ stages).

Given $\es \subseteq {\mathcal P}_{i}$ and a candidate branch set~$B$, we let
${\rm Eff}_\es(B)$ denote the set of all those effective vertices in~$B$ that are
matched to some (effective) vertex on a path in~$\es$ (in our random matching~$M^*$).
Assume that $|\es|= \alpha |{\mathcal P}_{i}|$ where $1/2\le \alpha\le 1$.
Let
\begin{equation}\label{eq:Ialpha}
I(\alpha):=\left[{\alpha t\over 3^i}\left(1\pm {1\over 4} \right)\right]
\end{equation}
and let~$\E_\es$ denote the event that $|{\rm Eff}_\es(B)|\in I(\alpha)$ for every candidate
branch set~$B$. Let $\overline{\E_\es}$ denote the complement of~$\E_\es$.

\begin{lemma}\label{Ei}
$\prob_i (\overline{\E_\es}) = \exp \left(-\Omega(\ln^2 n)\right)$.
\end{lemma}
\begin{proof}
Consider any candidate branch set~$B$. Note that $|{\rm Eff}_\es(B)|$ is
hypergeometrically distributed with mean $\lambda:= {\rm eff}_i(B)|\es|\ell_i/(kt-x_{i-1})$.
But $|\es|\ell_i=\alpha |\Pa_i| \ell_i = \alpha |Q_i|_{\rm eff}$ and
$kt-x_{i-1}={kt \over 2}(1+3^{-(i-1)})$ by~(\ref{eq:available}).
So~(\ref{eq:e_2j}) implies that
$$
{{\rm eff}_i(B)\over kt-x_{i-1}} \in \left[ {{t\over 2}(1+3^{-(i-1)}) \pm n^{1/3} \over
{kt\over 2}(1+3^{-(i-1)})}\right]
\subseteq \left[{1\over k}\left(1\pm {1\over 5} \right)\right]$$
and thus
$$ \lambda \in\left[ \alpha |Q_i|_{\rm eff}
\left({1\over k}\left(1\pm {1\over 5} \right)\right) \right] \stackrel{(\ref{eq:Qi})}{=}
\left[ {\alpha t\over 3^i}\left(1\pm {1\over 5} \right) \right].$$
In particular, together with (\ref{eq:Hyp_Conc1}) this implies that
\begin{equation*}
\prob_i \left(\big{|}\, |{\rm Eff}_\es(B)|-\lambda\big{|}
\geq n^{1/4} \ln n \right) =\exp\left(-\Omega(\ln^2 n)\right).
\end{equation*}
So with probability at most
$k \exp\left(-\Omega(\ln^2 n)\right) = \exp \left(-\Omega(\ln^2 n)\right)$
there is a candidate branch set~$B$ with $$|{\rm Eff}_\es(B)| \notin
\left[{\alpha t\over 3^i}\left(1\pm {1\over 5} \right) \pm n^{1/4}\ln n\right]
\stackrel{(\ref{eq:defkt}),(\ref{eq:Ialpha})}{\subseteq} I(\alpha ),$$
as required.
\end{proof}

\subsection{A lower bound for the degrees of the pairs of candidate branch sets in~$\cG_i$.}
Recall that, as described in the paragraph after~(\ref{eq:Bi}),
when analyzing the $i$th stage, we may assume that $|\U^*_{i-1}|< 26 U_{i-1}/27$
and thus $\U'_{i-1}$ is well defined.
Given a candidate branch set~$B$ and path~$P\in\Pa_i$, we write $P\sim B$
if some effective vertex on~$P$ is matched to some vertex in~${\rm Eff}_i(B)$ (in our
random matching~$M^*$). Consider an auxiliary bipartite
graph~$\Aux_i$ whose vertex classes are~$\U'_{i-1}$ and~$\Pa_i$ and in which a
pair $(B,B')\in \U'_{i-1}$ is adjacent to $P\in \Pa_i$ if~$P$ can be used to join $(B,B')$,
i.e.~if $P\sim B$ and $P\sim B'$. We will now estimate the
degrees of the vertices in ${\mathcal U}'_{i-1}$ in~$\cG_i$. Given $\es \subseteq {\mathcal P}_{i}$,
we let $d_{\cG_i}(L,\es)$ denote the degree of a vertex/pair $L \in {\mathcal U}'_{i-1}$ into
the set~$\es$ (in~$\cG_i$).

\begin{lemma}\label{lem:degL}
Suppose that $1/2 \le \alpha \le 1$ (where $\alpha$ may depend on $n$).
Fix $L \in {\mathcal U}'_{i-1}$ and $\es\subseteq \Pa_i$ with $|\es|= \alpha |{\mathcal P}_{i}|$. Then
$\prob_i \left(d_{\cG_i}(L,\es)\leq 1/(2\eps^{3}) \right)< 3\eps$.
\end{lemma}
\begin{proof}
Let $L=(B,B')$.
Our aim is to show that
\begin{equation}\label{eq:conditional}
\prob_i \left(d_{\cG_i}(L,\es)\leq 1/(2\eps^{3}) \ | \ \E_\es \right) < 2\eps.
\end{equation}
This implies the lemma since
\begin{eqnarray*}
\prob_i \left(d_{\cG_i}(L,\es)\leq 1/(2\eps^{3}) \right)
&\leq& \prob_i \left(d_{\cG_i}(L,\es)\leq 1/(2\eps^{3}) \ | \ \E_\es \right) + \prob_i (\overline{\E_\es}) \\
&\stackrel{(\ref{eq:conditional}),\ \mathtt{Lemma}~\ref{Ei}}{\leq} & 2\eps + \exp \left(-\Omega(\ln^2 n)\right)
< 3\eps.
\end{eqnarray*}
To estimate the number of all those paths in $\es$
that are neighbours of both $B$ and $B'$ in the auxiliary graph $\aux_i$,
we will first bound the number of paths in $\es$ that are neighbours of $B$ and
then we will estimate how many of them are neighbours of~$B'$.
More precisely, we
will first show that most of the paths $P\in \es$ with $P\sim B$
are joined to~$B$ by exactly one (matching) edge. Let us condition first on a particular
realization~$E_B$ of ${\rm Eff}_\es(B)$ with $|E_B|\in I(\alpha )$.
Denote the corresponding probability subspace of~$\prob_i$ (where we condition on
the event that ${\rm Eff}_\es(B)=E_B$ and on~$\E_\es$) by $\prob_{i,\E_\es,E_B}$.
Assuming an arbitrary ordering of the vertices in~$E_B$, we expose
their neighbours (in the random matching) on the paths in~$\es$ one by one according
to this ordering. We say that the $j$th vertex {\em fails} if its neighbour lies in a path from~$\es$
that already contains a neighbour of the previously exposed
vertices. Note that the number of paths in~$\es$ containing more than one
neighbour of~$E_B$ is bounded above by the number of failures that occur
during the exposure of the neighbours of~$E_B$. Suppose we have exposed the neighbours of the first
$j-1$ vertices in~$E_B$. Let the corresponding event be $\C_{j-1}$.
To estimate the probability that the $j$th vertex fails, observe
that the number of all those paths in~$\es$ that already have a neighbour in~$E_B$
is less than~$j$ and each of them contains less than~$\ell_i$ effective vertices
which are still available. Note that this holds regardless of
what~$\C_{j-1}$ is. Thus%
     \COMMENT{To see the last inequality it suffices to show that
$|E_B| \le \alpha |\Pa_i| \ell_i/2$. This holds if $|E_B| \ll |\Pa_i|$.
But $i \leq {1\over 6}\log_3 n$, so $|\Pa_i|=\Omega (n^{2/3})\gg t = \Theta(\sqrt{n})$.
Now note $|E_B| \le t$.}
\begin{align*}
\prob_{i,\E_\es,E_B}(\mbox{the $j$th vertex fails} \mid \C_{j-1}) <
{j\ell_i \over \alpha |\Pa_i|\ell_i - (j-1)}
\leq {|E_B|\ell_i \over \alpha |\Pa_i|\ell_i - |E_B|}
\stackrel{(\ref{eq:i0def}),(\ref{eq:U_isize})}{\leq} {2|E_B|\over \alpha |\Pa_i|} =:\tilde{p}.
\end{align*}
In particular, let $\D_{j-1}$ be any event which depends only on the neighbours
of the first $j-1$ vertices in~$E_B$. Then
\begin{equation} \label{ptilde}
\prob_{i,\E_\es,E_B}(\mbox{the $j$th vertex fails} \mid \D_{j-1}) \le \tilde{p}.
\end{equation}
Now let $A=\{a_1,\dots,a_r\}$ be any set of vertices in~$E_B$
(where~$a_q$ precedes~$a_{q+1}$ in the ordering of~$E_B$) and let Fail$_A$ denote the
event that the set of failure vertices equals~$A$. Then
$$
\prob_{i,\E_\es,E_B}(\mbox{Fail}_A)
\le \prod_{q=1}^r \prob_{i,\E_\es,E_B}( a_q \mbox{ fails } \mid a_1, \dots, a_{q-1} \mbox{ fail})
\stackrel{(\ref{ptilde})}{\le} \tilde{p}^r.
$$
This in turn implies that
\begin{align*}
\prob_{i,\E_\es,E_B}( \ge f \mbox{ failures})
 \le \sum_{r= f}^{|E_B|} \sum_{A \subseteq E_B \atop |A|=r}
\prob_{i,\E_\es,E_B}(\mbox{Fail}_A)
\le \sum_{r= f}^{|E_B|} \binom{|E_B|}{r} \tilde{p}^r
 \le \sum_{r= f}^{|E_B|} \left( \frac{e |E_B|}{r} \frac{2|E_B|}{\alpha|\Pa_i|} \right)^r.
\end{align*}
Since $|E_B|\in I(\alpha)$ we have
\begin{equation} \label{Tsquare}
\frac{|E_B|^2}{\alpha |\Pa_i|}\stackrel{(\ref{eq:U_isize}),(\ref{eq:Ialpha})}{\in}
\left[(1 \pm 1/4)^2 \frac{\alpha^2 t^2}{9^i} \frac{1}{\alpha}\frac{300\cdot 9^i}{9kt}\right]
\stackrel{(\ref{eq:defkt})}{\subseteq} \left[(1 \pm 2/3) \frac{100 \alpha}{3\sqrt{2} \eps^3}\right].
\end{equation}
Let $\F_B$ denote the event that at least $1000/\eps^3$ failures occur when we expose the neighbours
of the vertices in~$E_B$.
Thus by setting  $f := 1000/\eps^3$, we obtain
\begin{equation} \label{eq:ActualProbFault}
\prob_{i,\E_\es,E_B}(\F_B)=
\prob_{i,\E_\es,E_B}( \ge f \mbox{ failures})
\le \sum_{r \ge f} \left( \frac{5 \sqrt{2}\cdot 100 e\alpha}{9\eps^3r} \right)^r
\leq \sum_{r\geq f} (1/2)^r  \leq \eps.
\end{equation}
Let $\overline{\F_B}$ denote the complement of~$\F_B$. Note that if
$\overline{\F_B}$ occurs, then there are at least $|E_B|-1000 /\eps^3$
paths in~$\es$ that are joined to $E_B \subseteq B$ by exactly one (matching) edge. Let $\es (B)$
denote the set of these paths. So
\begin{equation}\label{MiBounds}
|E_B| - 1000 /\eps^3 \leq |\es (B)|\leq |E_B|
\end{equation}
(for the second inequality we need not assume that $\overline{\F_B}$ holds).
Now we additionally condition on a specific realization $E_{B'}$ of ${\rm Eff}_\es(B')$ with $|E_{B'}|\in I(\alpha)$.
As above, we fix an arbitrary ordering on the
vertices in~$E_{B'}$ according to which we expose their neighbours in~$\es$.
We say that the $j$th vertex of~$E_{B'}$ is {\em useful} if it is
adjacent to a vertex lying on a path from~$\es (B)$ such that none of
the previous vertices in~$E_{B'}$ is joined to this path. Note that if~$U(B')$ denotes the
set of vertices in $E_{B'}$ that are useful, then
\begin{equation} \label{eq:DegBound} |U(B')| \leq d_{\cG_i}(L,\es).
\end{equation}
Given $\overline{\F_B}$, we will show that with high probability $|U(B')|\geq  1/(2\eps^{3})$.

Note that there are exactly $R:=\alpha |\Pa_i|\ell_i - |E_B|$ effective vertices on the
paths in~$\es$ that are still available to be matched to the vertices of $E_{B'}$.
Put $s:=|\es (B)|(\ell_i-1)$ and let~$C$ be any subset of~$E_{B'}$ with $c:=|C| \leq  1/(2\eps^{3})$.
Suppose that~$C$ is the set of useful vertices in~$E_{B'}$. Then the vertices in~$C$ are
matched to effective vertices on different paths in~$\es(B)$. So there are $|\es (B)|_c (\ell_i-1)^c\le s^c$
choices for the neighbours of~$C$. Moreover, each vertex $x\in E_{B'}\setminus C$ is
either matched to an effective vertex on a path in~$\es(B)$ which already contains a neighbour
of~$C$ or~$x$ is matched to an effective vertex on a path in $\es\setminus \es(B)$.
There are less than $c\ell_i$ choices for a neighbour of~$x$ having the first property
and $R-s$ choices for a neighbour of~$x$ having the second property. Thus in total the
number of choices for the neighbours of $E_{B'}\setminus C$ is at most
\begin{align*}
\sum_{q=0}^{|E_{B'}|-c} (|E_{B'}|-c)_q(c\ell_i)^q(R-s)_{|E_{B'}|-c-q}
& \le (R-s)_{|E_{B'}|-c}\sum_{q=0}^{|E_{B'}|-c}\left(\frac{c\ell_i |E_{B'}|}{R/2}\right)^q\\
& \le (R-s)_{|E_{B'}|-c}\sum_{q\ge 0}\left(\frac{1}{2}\right)^q=2 (R-s)_{|E_{B'}|-c}.
\end{align*}
(Here we used that $|E_{B'}|=O(\sqrt{n})=o(|\Pa_i|)$ and so
$s=o(R)$ as well as $|E_{B'}|\ell_i=o(R)$.)
Setting $p:=s/R$ we obtain
\begin{align*}
\prob_{i,\E_\es,E_B} & (U(B') =C\ | \ \overline{\F_B}, {\rm Eff}_\es(B')=E_{B'} )
\le \frac{2s^c  \left(R-s \right)_{|E_{B'}|-c}}{ (R)_{|E_{B'}|}}\\
& \leq  2s^c \left( {R - s \over R} \right)^{|E_{B'}|-c}~\left( 1\over R-|E_{B'}| \right)^c
= 2\left({s\over R-|E_{B'}|}\right)^c \left(1- {s\over R} \right)^{|E_{B'}|-c}\\
& \le 2\left(\frac{10}{9}\right)^c\left(\frac{s}{R}\right)^c\left(1-\frac{s}{R}\right)^{|E_{B'}|-c}
=2 \left(\frac{10}{9}\right)^c p^c(1-p)^{|E_{B'}|-c}.
\end{align*}
(In the second inequality we used that ${a - j \over b - j}< {a \over b}$, for $0<j<a<b$ and
in the last inequality we again used that $|E_{B'}|=o(R)$.) Thus
\begin{align} \label{eq:useful}
\prob_{i,\E_\es,E_B} (|U(B')| \le 1/(2\eps^3) \mid & \ \overline{\F_B}, {\rm Eff}_\es(B')=E_{B'})
\nonumber\\ &\leq 2(10/9)^{1/(2\eps^{3})} \sum_{c \leq 1/(2\eps^{3})}
{|E_{B'}|\choose c} p^c (1-p)^{|E_{B'}|-c}.
\end{align}
Observe that the sum on the right-hand side is the probability that a binomial random
variable $Y$ with parameters $|E_{B'}|,p$
is at most $1/(2\eps^{3})$. To bound this probability from above, we will use the following
Chernoff bound (see e.g.~Inequality~(2.9) in \cite{JLR}):
\begin{equation} \label{chernoff}
\prob( Y \le \ex Y/2 ) \le 2\exp \left(-\ex Y/12 \right).
\end{equation}
Note that by (\ref{MiBounds}) and the definition of $R$, we have
$$p={s \over R}={|\es (B)|(\ell_i-1) \over R}
\geq {\left( |E_B| - 1000 /\eps^3\right) (\ell_i-1)
\over \alpha |\Pa_i|\ell_i} > { |E_B|\over \alpha|\Pa_i|}\frac{98}{100},
$$
where the last inequality holds since $\ell_i \ge 100$. Moreover,
the bound (\ref{Tsquare}) also holds if we replace
$|E_B|^2$ by $|E_B||E_{B'}|$. So altogether we have
\begin{eqnarray*}
\ex Y=|E_{B'}|p  \ge \frac{98}{100} \frac{|E_{B'}||E_B|}{\alpha |\Pa_i|}
 \stackrel{(\ref{Tsquare})}{\ge}
\frac{2\alpha}{\eps^3} \ge \frac{1}{\eps^3}.
\end{eqnarray*}
Thus~(\ref{chernoff}) implies that
\begin{equation*}
\sum_{c\leq 1/(2\eps^3)} {|E_{B'}|\choose c} p^c (1-p)^{|E_{B'}|-c} \leq
2\exp \left( - 1/(12\eps^3)\right).
\end{equation*}
Substituting this bound into (\ref{eq:useful}), we obtain%
    \COMMENT{Note that $(10/9)^7\le 2.1$}
$$
\prob_{i,\E_\es,E_B} ( |U(B')|\leq 1/(2\eps^3) \ | \ \overline{\F_B}, {\rm Eff}_\es(B')=E_{B'} )
\leq 4\left( (10/9)^{6}/e \right)^{1/(12\eps^3)} \le (9/10)^{1/(12\eps^3)}\le \eps.
$$
Since $E_{B'}$ was an arbitrary realization of~${\rm Eff}_\es(B')$ with
$|E_{B'}|\in I(\alpha)$, this implies that
$$
\prob_{i,\E_\es,E_B} (|U(B')| \le 1/(2\eps^{3}) \mid \overline{\F_B}) \le \eps
$$
and thus
\begin{eqnarray} \label{eq:ProbUsef}
\prob_{i,\E_\es,E_B} (|U(B')| \le 1/(2\eps^{3}))
& \le &
\prob_{i,\E_\es,E_B} (|U(B')| \le 1/(2\eps^{3}) \mid \overline{\F_B})
+\prob_{i,\E_\es,E_B} (\F_B)\nonumber\\
& \stackrel{(\ref{eq:ActualProbFault})}{\le} & 2\eps.
\end{eqnarray}
Finally, since~$E_B$ was an arbitrary realization of~${\rm Eff}_\es(B)$ with $|E_{B}|\in I(\alpha)$
it follows that
$$
\prob_i (d_{\cG_i}(L,\es)\le 1/(2\eps^{3}) \mid \E_\es)
 \stackrel{(\ref{eq:DegBound})}{\leq}  \prob_i (|U(B')|\le 1/(2\eps^{3}) \mid \E_\es)
\stackrel{(\ref{eq:ProbUsef})}{\leq}  2\eps,
$$
as required.
\end{proof}

Now given $\es \subseteq \Pa_i$, we let $\U(\es)$ be the set of all those pairs in~$\U'_{i-1}$
that have degree at most $1 /(2 \eps^3)$ into~$\es$ (in our auxiliary graph~$\cG_i$).
So if $1/2\le \alpha \le 1$ and $|\es|=\alpha |\Pa_i|$ then Lemma~\ref{lem:degL} implies that
\begin{equation} \label{T_expectation}
\ex_i (|\U(\es)|) \leq 3\eps|\U'_{i-1}|\le {U_{i-1}\over 2 \cdot 27},
\end{equation}
where $\ex_i (\cdot)$ denotes the expectation that arises from the probability measure
$\prob_i(\cdot)$.

\begin{lemma}\label{TBound}
Let $1/2\le \alpha \le 1$. Then every $\es\subseteq \Pa_i$ with $|\es|=\alpha |\Pa_i|$
satisfies
\begin{equation} \label{TBound1}
\prob_i \left(|\U(\es)| > {U_{i-1}\over 27} \right)
\leq  2\exp \left(- {2\eps^8 n \over (3^{i-1})^7} \right)
\end{equation}
as well as
\begin{equation} \label{Tbound2}
\prob_i \left(|\U(\es)| > {U_{i-1}\over 27} \right)
\leq  2\exp \left( -\eps^4 3^{(i-1)/4} \right).
\end{equation}
\end{lemma}
\begin{proof}
Our aim is to apply~(\ref{ConcM}) to show that~$|\U(\es)|$ is concentrated around its expected value.
We first prove~(\ref{TBound1}). Here~$W$ will be the space of all those matchings which match
the set~${\rm Eff}(Q_i)$ of effective vertices
on~$Q_i$ into the set $\bigcup_{j=1}^k {\rm Eff}_{i}(B_j)$ of all those effective vertices
in the candidate branch sets that
are still unmatched after the first~$i-1$ stages (equipped with the uniform distribution).
(Recall that ${\rm Eff}_{i}(B)$ is fixed since we condition on~$\M_{i-1}$.)
So each matching in~~$W$ consists of~$|Q_i|_{\rm eff}$ edges.
The metric~$d$ on~$W$ is defined by $d(M,M'):=2\ell_i|M\triangle M'|$
for all $M,M'\in W$. It is easy to see that this is indeed a metric.

So let us now define the partitions $F_0,\ldots,F_{|Q_i|_{\rm eff}}$. $F_0:=\{W\}$ and
each part of~$F_{|Q_i|_{\rm eff}}$ will consist of a single matching in~$W$. To define~$F_j$
for $1\leq j < |Q_i|_{\rm eff}$, fix a linear ordering on the vertices in~${\rm Eff}(Q_i)$.
Given a matching $M\in W$, the \emph{$j$-prefix of~$M$} is the set of all edges in~$M$
adjacent to the first~$j$ vertices in~${\rm Eff}(Q_i)$.
Each part of the partition~$F_j$ will consist of all those matchings in~$W$ having the
same~$j$-prefix. Clearly~$F_{j+1}$ refines~$F_j$.

To define the bijection~$\phi$, consider any two parts~$A\neq B$
of~$F_{j+1}$ and any part~$C$ of~$F_j$ such that $A,B\subseteq C$.
So if $M \in A$ and $M' \in B$, then~$M$ and~$M'$ have the same $j$-prefix and
they differ at the edge that is adjacent to the $(j+1)$th vertex in~${\rm Eff}(Q_i)$.
Let~$v_A$ and~$v_B$ be the neighbours of the $(j+1)$th vertex in~$M$ and~$M'$, respectively.
Note that~$v_A$ does not depend on the choice of $M\in A$ and similarly for~$v_B$.
We define $\phi: A \rightarrow B$ by saying that for all $M\in A$ the matching~$\phi(M)$
is obtained from~$M$ as follows: the $(j+1)$th vertex in~${\rm Eff}(Q_i)$ is now matched to~$v_B$
and~$v_A$ is matched to the neighbour of~$v_B$ in~$M$, every other edge of~$M$ remains
unchanged. Thus the size of the symmetric difference of~$M$ and~$\phi (M)$ is 4
and so $d(M, \phi (M))\leq 8\ell_i$. So we can take $c:=8\ell_i$.

Now note that~$|\U(\es)|$ is a function whose value is determined by a matching from~$W$
chosen uniformly at random. So we take~$f: W\rightarrow\mathbb{R}$ to be the function
defined by setting $f(M)$ to be the value of~$|\U(\es)|$ on~$M$ (for all $M\in W$).
We have to show that for any $M,M' \in W$ we have $|f(M)-f(M')|\leq d(M,M')$.
To do so, we will construct a sequence $M_0,M_1,\ldots, M_{q}$ of matchings in~$W$ such that
$M_0:=M$, $M_q:=M'$ and such that~$M_j$ and~$M'$ agree on the first~$j$ vertices in~${\rm Eff}(Q_i)$
(i.e.~$M_j$ and $M'$ have the same prefix).
Suppose that we have constructed~$M_j$ for some~$j<q$ and that we
now wish to construct~$M_{j+1}$. Let~$v$ be the first vertex in~${\rm Eff}(Q_i)$ on
which~$M_j$ and~$M'$ differ. Let~$b$ be its neighbour in~$M'$ and let~$v'$
be the neighbour of~$b$ in~$M_j$. Define~$M_{j+1}$ to be the matching obtained from~$M_j$ by
swapping the neighbours of~$v$ and~$v'$ in~$M_j$. So~$M_{j+1}$ now agrees with~$M'$
on~$v$ and all (the at least~$j$) vertices preceding~$v$ in~${\rm Eff}(Q_i)$.
Note that $|f(M_j)-f(M_{j+1})|\leq 4\ell_i$ since swapping two edges can change~$|\U(\es)|$
by at most $4\ell_i$. Indeed, to see the latter, note that for each one of these two
edges there are $\ell_i - 1$ other edges starting from the same path in~$\Pa_i$,
and therefore each of these two edges contributes to the degree of at most~$\ell_i$ pairs
in~$\U'_{i-1}$. If we swap these edges, this might change the degree
of at most~$4\ell_i$ pairs. So~$|\U(\es)|$ can be increased or decreased by
at most~$4\ell_i$. Also observe that $q\le |M \triangle M'|/2$
since initially the number of vertices in~${\rm Eff}(Q_i)$ on which~$M$ and~$M'$
differ equals $|M \triangle M'|/2$ and in each step this number decreases by at least~1.
Therefore,
\begin{equation} \label{eq:bbd_diff}
|f(M)-f(M')|\leq \sum_{j=0}^{q-1} |f(M_{j})-f(M_{j+1})| \leq 4q\ell_i
 \leq 2\ell_i |M\triangle M'|= d(M,M').
\end{equation}
Now, we are ready to apply~(\ref{ConcM}): if $|\U(\es)| > U_{i-1}/27$, then by~(\ref{T_expectation}) we have
$|\U(\es)|- \ex (|\U(\es)|)>U_{i-1}/(2\cdot 27)$ and~(\ref{ConcM}) yields
\begin{eqnarray*}
\prob_i \left(|\U(\es)| > {U_{i-1}\over 27} \right) & \leq &
 2\exp \left(-2 {(U_{i-1}/(2\cdot 27))^2\over |Q_i|_{\rm eff}8^2\ell_i^2} \right) \\
& \stackrel{(\ref{eq:Qi}),(\ref{eq:lidef}),(\ref{eq:InductiveStep})}{=} &
  2 \exp \left( - \frac{1}{2^7 \cdot 27^2} \frac{\eps^8 n^2}{(27^{i-1})^2}
\frac{3^i}{(1+o(1))\sqrt{2}\eps n} \frac{1}{100^2 9^{i-1}} \right) \\
& \le & 2 \exp \left( -\frac{2\eps^8 n}{(3^{i-1})^7} \right).
\end{eqnarray*}
Now we prove~(\ref{Tbound2}). In this case we can apply~(\ref{ConcM}) with metric
$d(M,M'):=2\Delta_{i-1}|M\triangle M'|$ and $c:=8\Delta_{i-1}$. Indeed,%
     \COMMENT{Here swapping 2 edges will change the degree of at most $2\Delta_{i-1}$
pairs. So $d(M,M'):=\Delta_{i-1}|M\triangle M'|$ and $c:=4\Delta_{i-1}$ would also
work. But perhaps it is clearer to make things similar to the 1st case.}
for each candidate branch set~$B$
and each $P\in\Pa_i$ the removal/addition of a $B$-$P$ edge can only affect the degrees of those
pairs in~$\U'_{i-1}$ which contain~$B$. But there are at most $\Delta_{i-1}$ such pairs. Thus
\begin{eqnarray*}
\prob_i \left(|\U(\es)| > {U_{i-1}\over 27} \right) & \leq &
 2\exp \left(-2 {(U_{i-1}/(2\cdot 27))^2\over |Q_i|_{\rm eff} 8^2\Delta_{i-1}^2} \right) \\
 & \stackrel{(\ref{eq:Qi}),(\ref{eq:def_of_delta})}{\le} &2\exp \left(-
\frac{ U_{i-1}^2}{2^7\cdot 27^2} \frac{3^i}{kt}
\frac{\eps^{1/4} k^2}{ (3/2)^{2(i-1)} U_{i-1}^2} \right) \\
& \stackrel{(\ref{eq:defkt})}{\le} & 2\exp \left( -\frac{1}{2^7 \cdot 27^2}
\frac{3\eps^{1/4} (1+o(1))\sqrt{2}\eps^2}{1/\eps} (4/3)^{i-1} \right). \\
& \le & 2\exp \left( - \eps^4 (4/3)^{i-1} \right) \le 2\exp \left( - \eps^4 3^{(i-1)/4} \right),
\end{eqnarray*}
as required.
\end{proof}

Define $\beta$ by
\begin{equation}\label{eq:defa0}
\beta:= \left( \eps^{8}\over 2 (3^{i-1})^{7}   \right)^2.
\end{equation}

\begin{lemma} \label{bad_subsets}
For each $i$ with $3^{i-1} \le n^{2/33}$ let~$Y_i$ denote the number of all
those subsets~$\es$ of~$\Pa_i$ with $|\es|=(1-\beta) |\Pa_i|$ for which $|\U(\es)|>U_{i-1}/27$.
Then $\prob_i(Y_i>0)\le n^{-1/34}$.
\end{lemma}
\begin{proof}
Note that~(\ref{eq:U_isize}) and the restriction on~$i$ together imply that
$\beta|\Pa_i|=\Omega(n/3^{16i})=\Omega(n^{1/33})$ and so we may treat it as an integer.
(\ref{TBound1}) implies that
\begin{equation} \label{eq:exp}
\ex_i(Y_i)\le {|\Pa_i|\choose (1-\beta) |\Pa_i|}
2\exp \left(- {2\eps^8 n \over (3^{i-1})^7} \right).
\end{equation}
Note that $|\Pa_i| \le n$. So
$$
{|\Pa_i|\choose (1-\beta) |\Pa_i|}={|\Pa_i|\choose \beta |\Pa_i|}
\leq  \left( {e \over \beta} \right)^{\beta |\Pa_i|}
\leq \beta^{-2\beta|\Pa_i|} \le \beta^{-2\beta n}.
$$
Now note that if $a>0$ is sufficiently small then%
     \COMMENT{Note that $x\ln x^{-1}/x^{1/2} \rightarrow 0$ as $x \rightarrow 0$.
This is the case since $x^{1/2} \ln x^{-1}\rightarrow 0$
as $x \rightarrow 0$. Indeed $x^{1/2}\ln x^{-1}=\exp \left({1\over 2}\ln x +\ln \ln x^{-1}\right)=
\exp \left(-{1\over 2}\ln x^{-1} +\ln \ln x^{-1}\right) \rightarrow 0$, as $x \rightarrow 0$.}
$a\ln (a^{-1})\leq a^{1/2}$. Thus
$$
{|\Pa_i|\choose (1-\beta)|\Pa_i|}\le e^{2\beta^{1/2}n} \stackrel{(\ref{eq:defa0})}{=}
\exp \left(  \frac{\eps^8 n}{(3^{i-1})^7} \right).
$$
So if $3^{i-1} \le n^{2/33}$ then
\begin{eqnarray*}
\prob_i(Y_i>0) \le \ex_i(Y_i)\stackrel{(\ref{eq:exp})}{\le} 2 \exp \left(-{\eps^8 n \over (3^{i-1})^7} \right)
= \exp\left(-\Omega \left( n^{19/33} \right)\right)\le n^{-1/34},
\end{eqnarray*}
as required.
\end{proof}


\subsection{An upper bound on the degrees of the paths in~$\cG_i$}
Let $d:=10^6$. We now estimate the probability that a given path $P \in \Pa_i$ joins at least~$d$
unjoined pairs in $\U'_{i-1}$.
\begin{lemma}\label{prob_big_deg}
If $d=10^6$, $i\le i_0$ and~$i$ satisfies
\begin{equation} \label{eq:loweri}
3^{i-1} \ge 1/\eps,
\end{equation}
then for every fixed $P\in \Pa_i$ we have $ \prob_i (d_{\cG_i}(P)\geq d) \leq \beta/3^i$.
\end{lemma}
\proof
Suppose that $C\subseteq \U'_{i-1}$ is a set of size~$d$ which lies in the neighbourhood
of~$P$ in the auxiliary graph~$\cG_i$. Let $\B(C)$ denote the set of candidate branch sets involved
in the pairs from~$C$. Note that
$$\sqrt{2d}\le |\B(C)|\le 2d.$$
Moreover, $P\sim B$ for each candidate branch set~$B\in \B(C)$. (Recall that this means that
there is an effective vertex on~$P$ that is matched to some vertex in~${\rm Eff}_i(B)$, where
${\rm Eff}_i(B)$ was the set of all those effective vertices in~$B$
that are still available after the $(i-1)$th stage.)
Now let~${\bf B}$ be the collection of all the sets~$\B$ of
candidate branch sets such that for each $B\in \B$ there is a $B'\in \B$ with
$(B,B')\in \U'_{i-1}$ and such that $b:=|\B|$ satisfies
\begin{equation}\label{eq:B}
\sqrt{2d}\le b\le 2d.
\end{equation}
Thus $\B(C)\in {\bf B}$ for any~$C$ as above and hence
\begin{equation}\label{eq:degP}
\prob_i (d_{\cG_i}(P)\geq d)\le \sum_{\B\in {\bf B}}\prob_i (P\sim B \ \forall B\in \B).
\end{equation}
To bound the latter probability, consider any $\B\in {\bf B}$, let $b:=|\B|$ and
$s:=\sum_{B\in \B} {\rm eff}_i(B)\leq b t$. Recall that~$x_{i-1}$ was the total number of
all those effective vertices in the branch sets that are matched after the first~$i-1$ stages.
So
\begin{eqnarray*}
\prob_i (P\sim B \ \forall B\in \B)
&\leq & {\ell_i \choose b} {(s)_{b}\over (kt-x_{i-1})_{b}}
\leq \left( {e \ell_i  \over b} \right)^{b} \left({ s\over kt- x_{i-1}}\right)^{b} \\
& \le &
\left( {e \ell_i  \over b} { b t\over kt- x_{i-1}}\right)^{b}
\stackrel{(\ref{eq:available})}{\leq}
\left( { e t \ell_i \over kt/2}\right)^{b}  = \left( {2e  \ell_i \over k} \right)^{b}.
\end{eqnarray*}
In the second inequality, we used that ${a - j \over b - j}< {a \over b}$, for $0<j<a<b$.
To bound~$|{\bf B}|$, consider an auxiliary graph $\AuxU_{i-1}$
whose vertex set is the set of candidate branch sets and whose edges correspond to the
pairs in~$\U_{i-1}'$. Since~$\AuxU_{i-1}$ involves only edges/pairs from~$\U_{i-1}'$
its maximum degree is at most~$\Delta_{i-1}$. Consider any~$b$ as in~(\ref{eq:B}).
Note that each $\B\in {\bf B}$ with~$|\B|=b$ corresponds to a subgraph~$F$ of~$\AuxU_{i-1}$
which has order~$b$ and in which no vertex is isolated. We claim that for all $q \le b/2$, the number of
such subgraphs~$F$ having precisely~$q$ components is at most $U_{i-1}^q (b\Delta_{i-1})^{b -2q}$.
To see this, note that each component of~$F$ has to contain at least one edge
(this is also the reason why it makes sense only to consider $q \le b/2$). So
each subgraph $F$ as above can be obtained as follows. First choose~$q$ (independent) edges of
$\AuxU_{i-1}$. The number of choices for this is at most $U_{i-1}^q$.
Now successively add the remaining $b-2q$ vertices to the existing subgraph without creating
new components. At each step there are at most~$b$ vertices~$y$ to which a new
vertex~$z$ can be attached and once we have chosen $y$, there are at most
$\Delta(\AuxU_{i-1})\le \Delta_{i-1}$ choices for~$z$, which proves the claim.

Let~${\bf B}_{b,q}$ be the set of all those $\B\in {\bf B}$ that have size~$b$
and  induce~$q$ components in~$\AuxU_{i-1}$. Then
\begin{eqnarray*}
 & & \sum_{\B\in {\bf B}_{b,q}} \prob_i(P\sim B \ \forall B\in \B)
 \stackrel{(\ref{eq:B})}{\leq}
U_{i-1}^q \left(2d \Delta_{i-1}\right)^{b-2q} \left( {2e\ell_i \over k} \right)^{b}\\
& \stackrel{(\ref{eq:lidef}), (\ref{eq:InductiveStep}), (\ref{eq:def_of_delta})}{\leq} &
\left( \frac{U_0}{27^{i-1}} \frac{1}{4d^2} \frac{\eps^{1/4} k^2 81^{i-1}4^{i-1}}{U_0^2} \right)^q
\left( 2d \frac{U_0 }{\eps^{1/8} k 9^{i-1}2^{i-1}} \frac{2e \cdot 100 \cdot 3^{i-1}}{k} \right)^{b}\\
& \stackrel{(\ref{eq:U_0Value})}{\le } &
(12^{i-1})^q \left( \frac{2000d}{ \eps^{1/8}6^{i-1}} \right)^{b}
\le \left( \frac{4 \cdot 10^6 d^2}{ \eps^{1/4} 3^{i-1}} \right)^{b/2}
\le \left( \frac{1}{ \eps^{1/3} 3^{i-1}} \right)^{b/2}
\stackrel{(\ref{eq:loweri})}{\le} \left( \frac{ \eps^{1/6}}{ 3^{(i-1)/2}} \right)^{b/2}.
\end{eqnarray*}
Since $b \geq \sqrt{2 d} \ge 4 \cdot 48$ by~(\ref{eq:B}) this implies
\begin{align*}
\prob_i (d_{\cG_i}(P)\ge d) &\stackrel{(\ref{eq:degP})}{\le}
\sum_{\sqrt{2d}\le b\le 2d}\ \sum_{q=1}^{b/2}\left( \frac{ \eps^{1/6}}{ 3^{(i-1)/2}} \right)^{b/2}
\le 2d^2 \left( \frac{ \eps^{1/3} }{ 3^{i-1}}\right)^{\sqrt{2d}/4}
\le 2d^2 \left( \frac{\eps^{1/3}}{ 3^{i-1}} \right)^{48}\le {\beta \over 3^i},
\end{align*}
as required.
\endproof

Given $d=10^6$ and $i$ satisfying~(\ref{eq:loweri}),
let~$\es_{i}$ denote the set of paths in~$\Pa_i$ which have degree less than~$d$ in~$\cG_i$.
We will now use Lemma~\ref{prob_big_deg} to show that with high probability~$\es_{i}$ is large.

\begin{lemma} \label{Exp_Onthe_Right}
$\prob_i \left( |\es_i|\geq (1-\beta) |\Pa_i | \right) \ge 1-n^{-1/34}$
for every $i\leq i_0$ which satisfies~(\ref{eq:loweri}).
\end{lemma}
\begin{proof}
Let $\bar{\es_i}:=\Pa_i \setminus \es_i$.
Note that Lemma~\ref{prob_big_deg} implies
\begin{equation*} 
\ex_i (\bar{\es_i})  = \prob_i (d_{\cG_i}(P)\geq d) |\Pa_i|\le {\beta \over 3^i} |\Pa_i|.
\end{equation*}
If $ (\log_3 n)/34 \leq i \leq i_0$ then together with Markov's inequality this yields
\begin{eqnarray*}
\prob_i (|\bar{\es}_i| > \beta|\Pa_i| )\leq {1\over 3^{i}}\le {1\over n^{1/34}}
\end{eqnarray*}
and thus Lemma~\ref{Exp_Onthe_Right} holds for all such~$i$.

If $i\le (\log_3 n)/34$ we will use~(\ref{ConcM}).
As in the proof of Lemma~\ref{TBound}, the underlying metric space will
be the set of all those matchings which match the set~${\rm Eff}(Q_i)$ of effective vertices
on~$Q_i$ into the set $\bigcup_{j=1}^k {\rm Eff}_{i}(B_j)$ of all those effective vertices
in the candidate branch sets that are still unmatched after the first~$i-1$ stages.
The series of partitions and the bijections~$\phi$ are also as defined there.
However, the metric imposed on~$W$ now changes: for any $M,M' \in W$ we set $d(M,M')=|M\triangle M'|$.
In particular this means that we can take $c:=4$.
$f:W\to \mathbb{R}$ will be the function defined by taking~$f(M)$ to be the value of~$|\bar{\es}_i|$
on~$M$ (for all $M\in W$).
Note that the analogue of (\ref{eq:bbd_diff}) is satisfied, since if we switch the
endpoints of two edges of a matching (as it is the case when we
obtain~$M_{j+1}$ from~$M_j$ as in the proof of Lemma~\ref{TBound}) $|\bar{\es}_i|$ changes by at most~2
as switching two edges only affects the degree of the (at most) two paths involved. Thus
\begin{equation*}
|f(M)-f(M')|\leq \sum_{j=0}^{q-1} |f(M_{j})-f(M_{j+1})| \leq 2q
 \leq |M\triangle M'|\leq d(M,M').
\end{equation*}
Hence applying (\ref{ConcM}) with $a:= \beta|\Pa_i|/2$ we obtain
\begin{equation*}
\prob_i (|\bar{\es}_i|\geq \beta|\Pa_i|)
\leq 2 \exp \left(-2\frac{\beta^2 |\Pa_i|^2}{4\cdot 16 |Q_i|_{\rm eff}}\right).
\end{equation*}
So to complete the proof, it suffices to show that
$\beta^2 |\Pa_i|^2/|Q_i|_{\rm eff}= \Omega \left(n^{3/34} \right)$,
as this gives an error bound of $\exp ( -\Omega(n^{3/34}) ) \le 1/n^{34}$.
To prove the former, note that by (\ref{eq:Qi}), (\ref{eq:U_isize}) and
(\ref{eq:defa0}) we obtain
\begin{eqnarray*}
\frac{\beta^2 |\Pa_i|^2}{|Q_i|_{\rm eff}}
= \Theta \left( {1\over (3^{i-1})^{28}}~
\left({n \over 9^{i-1}} \right)^2~{3^i\over n}\right)
= \Theta \left({n \over 3^{31(i-1)} }\right)
= \Omega \left({n \over 3^{31(\log_3 n) /34} }\right)
= \Omega \left(n^{3/34} \right),
\end{eqnarray*}
as required.
\end{proof}


\subsection{Finding a large matching of~$\cG_i$}\label{sec:final}
The next lemma shows that with high probability we can join the required number of pairs from~$\U'_{i-1}$
during the $i$th stage.
\begin{lemma} \label{hall}
For each $i \le i_0$ we have $\prob_i (|\U_i| > U_{i-1}/27 ) \le 2n^{-1/34}$.
\end{lemma}
\proof
Recall that~$\U'_{i-1}$ (defined after~(\ref{eq:Bi})) was obtained from~$\U_{i-1}$ by
discarding all those pairs containing a candidate branch set from~$\B_i$. By definition
of~$\U_i$ we may assume that $|\U'_{i-1}|\ge U_{i-1}/27$ and it suffices to show
that in $\cG_i$ we can find a matching which covers all but at most $U_{i-1}/27$ vertices/pairs
in~$\U_{i-1}'$.

\medskip

\noindent
{\bf Case 1:} $3^{i-1} < 1/\eps$.\\
In this case, we apply~(\ref{TBound1}) with $\es:=\Pa_i$ (i.e.~$\alpha=1$) to obtain
that with probability at least
$1-2\exp \left(- {2\eps^8 n /(3^{i-1})^7} \right) \ge 1-2n^{-1/34}$ we have the
following: there is a set $\W \subseteq \U_{i-1}'$
with $|\W| = |\U_{i-1}'|-U_{i-1}/27$ so that every pair in~$\W$ has
degree at least $1/(2\eps^3)$ in $\cG_i$. On the other hand, clearly every path in~$\Pa_i$ has degree at most
$\ell_i^2 =10^4 9^{i-1}< 10^4/\eps^2$ in $\cG_i$. This implies that the subgraph
$\cG_i'$ of $\cG_i$ induced by $\W$ and $\Pa_i$ has a matching covering all of~$\W$.
To see this, consider any $\W' \subseteq \W$ and let $N(\W') \subseteq \Pa_i$ denote its neighbourhood
in $\cG_i'$. Then by counting edges between $\W'$ and $N(\W')$ we obtain that
$|\W'| /(2\eps^3) \le  (10^4/\eps^2) |N(\W')|.$
This in turn implies that $|N(\W')| \ge |\W'|$ and so Hall's condition is satisfied.
But this means that we can take $\U_i:=\U_{i-1}' \setminus \W$. Note that
$U_i=|\U_{i-1}'|-|\W|=U_{i-1}/27$, as required.
\medskip

\noindent
{\bf Case 2:} $1/\eps \le 3^{i-1} \le n^{2/33}$.\\
In this case we first apply
Lemma~\ref{Exp_Onthe_Right} to see  that with probability
at least $1- n^{-1/34}$ we have $|\es_i| \ge (1-\beta)|\Pa_i|$. By taking a subset we may
assume that $|\es_i|= (1-\beta)|\Pa_i|$. On the other hand, Lemma~\ref{bad_subsets} implies that
with probability at least $1-n^{-1/34}$ any set~$\es$ of this size satisfies $|\U(\es)|\le U_{i-1}/27$.
So with probability at least $1-2n^{-1/34}$ we have $|\U(\es_i)|\le U_{i-1}/27$. But if this
is the case then there is a set $\W \subseteq \U'_{i-1}$ with $|\W| = |\U'_{i-1}|-U_{i-1}/27$ so
that every pair in~$\W$ has degree at least $1/(2\eps^3)$ in the subgraph~$\cG_i''$
of~$\cG_i$ induced by~$\W$ and~$\es_i$.
On the other hand, the definition of~$\es_i$ implies that in~$\cG_i''$, the degree
of every vertex in~$\es_i$ is at most $d=10^6$.
As in the previous case, this implies that~$\cG_i''$ has a matching covering all of~$\W$. Indeed, to
verify Hall's condition consider any $\W' \subseteq \W$ and let $N(\W') \subseteq \es_i$
denote its neighbourhood in~$\cG_i''$. Then $|\W'| /(2\eps^3) \le  10^6 |N(\W')|$.
As before, we can take $\U_i:=\U_{i-1}' \setminus \W$.

\medskip

\noindent
{\bf Case 3:} $3^{i-1} \ge n^{2/33}$.\\
In this case, we apply~(\ref{Tbound2}) to $\es:=\Pa_i$ in order to obtain that with
probability at least
$$
1- 2\exp \left( -\eps^4 3^{(i-1)/4} \right)
\ge 1- 2\exp \left( -\eps^4 n^{1/66} \right) \ge 1-n^{-1/34}
$$
we have the following: there is a set $\W \subseteq \U_{i-1}'$
with $|\W| = |\U_{i-1}'|-U_{i-1}/27$ so that every pair in~$\W$ has
degree at least $1/(2\eps^3)$ in $\cG_i$. On the other hand, Lemma~\ref{prob_big_deg}
implies that the probability that $\Pa_i$ does not contain a path of degree at least $d=10^6$
in~$\cG_i$ is at least
$$1- \frac{\beta |\Pa_i|}{3^i} \stackrel{(\ref{eq:U_isize}),(\ref{eq:defa0})}{\ge}
 1- O(n/3^{17 i})=1- O(n^{-1/33})\ge 1-n^{-1/34}.
$$
So we may assume that both events occur and we get a matching covering all of~$\W$ in the subgraph~$\cG'_i$
of~$\cG_i$ induced by~$\W$ and~$\Pa_i$
as before. So we again obtain a set~$\U_i$ of the desired size, with the required error bounds.
\endproof

To complete the proof of Theorem~\ref{Thm_3Reg} it remains to combine all the error
probabilities for all the~$i_0=(\log_3 n)/6$ stages.
Recall that when analyzing the~$i$th stage we conditioned on~$\M^*_{i-1}$ (defined
after~(\ref{eq:e_2j})). However, all our probability bounds hold regardless of what the actual value
of~$M^*_{i-1}$ is (as long as~$M^*_{i-1}$ is good). So suppose that $|\U_{i-1}|=U_0/27^{i-1}$
for some $i\le i_0$. If $|\U_i|\neq U_0/27^i$ then either some candidate branch set
violated~(\ref{eq:e_2j}) or we had the undesired event that $|\U_i|>U_{i-1}/27$ in Lemma~\ref{hall}.
Thus
$$\prob(|\U_i|= U_0/27^i\mid |\U_{i-1}|=U_0/27^{i-1}, \X,\X')\ge 1-\exp(-\Omega (\ln^2 n))-2n^{-1/34}
\ge 1-3n^{-1/34}
$$
and so
$$\prob(|\U_i|= U_0/27^i \ \forall i\le i_0 \mid \X,\X')\ge 1-3i_0n^{-1/34}\ge 1-n^{-1/35}.
$$
This bound holds regardless of what the choices of $X_1,X_2,X'_1,X'_2$ actually are (as long as
$|X_1|=|X_2|$ is within the range determined in Lemma~\ref{Xrange}).
The only other reason why $|\U_{i_0}|\neq U_0/27^{i_0}$ is that we had an undesired event in
Lemma~\ref{Xrange}. This happens with probability $O(1/\ln ^2 n)$.
Altogether this shows that with probability $1-n^{-1/35}-O(1/\ln ^2 n)=1-o(1)$ after the $i_0$th stage
we are left with
\begin{equation} \label{eq:unjoined} U_{i_0} = {U_0\over 27^{i_0}}
\stackrel{(\ref{eq:i0def}),(\ref{eq:U_0Value})}{=} {\eps^4 n^{1/2}}
\end{equation}
unjoined pairs. We now discard a candidate branch set in each of these pairs
as well as all the candidate branch sets in $\B_0\cup\dots\cup \B_{i_0-1}$.
By~(\ref{eq:TotalB_i}) and~(\ref{eq:unjoined}) this gives a complete minor on
$k-6\eps^{1/8}k-{\eps^4 n^{1/2}} \ge \eps^2 n^{1/2}$
vertices, as required.

\section{Acknowledgements}
We are grateful to Tomasz {\L}uczak for helpful discussions on the phase transition of $G_{n,p}$.

\medskip

{\footnotesize \obeylines \parindent=0pt

Nikolaos Fountoulakis, Daniela K\"{u}hn \& Deryk Osthus
School of Mathematics
University of Birmingham
Edgbaston
Birmingham
B15 2TT
UK
}

{\footnotesize \parindent=0pt

\it{E-mail addresses}:
\tt{\{nikolaos,kuehn,osthus\}@maths.bham.ac.uk}}

\end{document}